\documentclass[12pt]{article}

\usepackage{amssymb}
\usepackage{amsmath}

\def\N{\mathbb N}

\def\E{\mathbf E}
\def\p{\mathbf P}
\def\d{{\rm d}}
\def\e{{\rm e}}

\newtheorem{theorem}{Theorem}

\newenvironment{proof*}{\smallskip \noindent {\bf Proof.} \ }
{\hspace*{1pt} \hfill\rule{6pt}{6pt} \medskip}

\newenvironment{proof}[1]{\smallskip \noindent {\bf Proof of #1.} \ }
{\hspace*{1pt} \hfill\rule{6pt}{6pt} \medskip}

\newcounter{example} \newcounter{refexample}
\newenvironment{example}{\smallskip \noindent \stepcounter{example} \refstepcounter{refexample}
{\bf Example \theexample.} \ } { 
\medskip}

\begin{document}

\begin{center}

\textbf{
\Large{
Asymptotics of nearly critical Galton--Watson processes with immigration}}

\vspace{1 cm}

\textsc{P\'eter Kevei
\footnote{The research was supported by the Analysis and Stochastics
Research Group of The Hungarian Academy of Sciences.}}

\vspace{1 cm}

Centro de Investigaci\'on en Matem\'aticas,
Jalisco S/N, Valenciana, Guanajuato, GTO 36240, Mexico;
\textit{e-mail:} \texttt{kevei@cimat.mx}

\end{center}

\bigskip

\begin{abstract}
We investigate the inhomogeneous Galton--Watson processes with immigration, where
$\rho_n$ the offspring means in the $n^\textrm{th}$ generation tends to 1. 
We show that if the second derivatives of the offspring generating functions go
to 0 rapidly enough, then the asymptotics are the same as in the INAR(1) case,
treated in \cite{GyIPV}. We also determine the limit if this assumption does not
hold showing the optimality of the conditions.

\end{abstract}

\section{Introduction}

Let $X_0= 0$ and consider the following inhomogeneous Galton--Watson process
with immigration:
$$
X_n = \sum_{j=1}^{X_{n-1}} \xi_{n,j} + \varepsilon_n,
$$
where $\{\xi_{n,j}, \varepsilon_n: n, j \in \N \}$ are independent nonnegative integer
valued random variables such that $\{\xi_{n,j} : j \in \N \}$ are identically distributed.
If the offspring distribution is Bernoulli distribution, that is each particle either
dies without descendant or leaves exactly one descendant, we obtain the so-called
first order integer-valued autoregressive (INAR(1)) time series.
We assume that the process is nearly critical, that is
$$
\E \xi_{n,1} = \rho_n \uparrow  1 \quad \textrm{as } n \to \infty.
$$
We note here that in the followings any non-specified limit relation is meant
as $n \to \infty$.

The theory of branching processes is used to model the evolution
of populations whose members live, reproduce and die independently of each
other. Its first appearance was motivated by the problem of extinction probability
of the family names in the British peerage. The problem was addressed by Francis Galton in
1873 and solved by Henry Watson in 1874. Since then the theory is developing to model
more complex systems, and now branching processes play an important role in models of genetics,
molecular biology, physics and computer science. Therefore we do not even
try to give a comprehensive bibliography. As a main reference on branching processes
we refer to the classical book of Athreya and Ney \cite{athreya-ney}.
For some recent application of INAR models and branching processes see \cite{GyIPV}. 

The aim of the present paper is to investigate the asymptotic properties of
nearly critical Galton--Watson processes with immigration
under general offspring distribution, and thus extend the results
in the INAR(1) case by Gy\"orfi, Isp\'any, Pap and Varga \cite{GyIPV},
by dropping the restrictive condition of Bernoulli offsprings.

Section 2 contains all the results and the discussions. In subsection 2.1 we investigate
the case which is parallel to the results in \cite{GyIPV}. Here we assume that the
variance of the offsprings tends to 0 with a prescribed rate, and thus in a clear
sense the offsprings are `almost' Bernoulli random variables. The methods here
are similar the ones in \cite{GyIPV}, however we emphasize that the proofs
are necessarily more difficult, since we do not have a closed form for the
generating functions even  in the case of Bernoulli immigration.
Subsection 2.2 deals with a significantly different case, when
the second derivatives also contribute to the limit distribution, and so the proofs need
new ideas. Subsection 2.3 contains  the case, when $\rho_n \to 1$ very fast.
Finally, in subsection 2.4 we consider the case of linear fractional generating
functions. This example is very important because of the explicit computations.
All the proofs are placed in section 3.

\section{Results and discussion}

Introduce the generating functions
$$
F_n(x) = \E x^{X_n}, \quad G_n(x) = \E x^{\xi_{n,1}}, \quad H_n(x) = \E x^{\varepsilon_n},
\ x \in [0,1].
$$
Using the branching property we  obtain the basic recursion for $F_n$:
\begin{eqnarray*}
F_n(x) & = & \E \left[ x^{\sum_{i=1}^{X_{n-1}} \xi_{n,i}  + \varepsilon_n }\right]
= \E \left[ \E \left( x^{\sum_{i=1}^{X_{n-1}} \xi_{n,i}  + \varepsilon_n }
\bigg | X_{n-1} \right) \right] \\
& = & \E \left[ G_n(x)^{X_{n-1}} \right] H_n(x) =
F_{n-1}(G_n(x)) H_n(x).
\end{eqnarray*}
Introduce the following notation:
$\overline G_{n+1,n}(x) = x$ and if $\overline G_{j+1,n}(x)$ is defined for $j \leq n$ then
$$
\overline G_{j,n}(x) =  \left( G_j \! \circ \ldots  \circ G_n \right)(x)
= G_j(\overline G_{j+1,n}(x)).
$$
With this notation the induction above gives the formula
\begin{equation} \label{gen-fgv}
F_n(x)= \prod_{j=1}^n H_j( \overline G_{j+1, n}(x) ).
\end{equation}

According to the continuity theorem for discrete random variables
(\cite{feller} p.~280) for proving a limit theorem we have to show
that $F_n(x)$ converges as $n \to \infty$ for all $x \in (0,1)$, and
the limit function is the generating function of the limit distribution.
Since for fix $j$ the function $H_j(\overline G_{j+1,n}(x)) \sim 1$,
we introduce the corresponding generating function
$$
\tilde F_n (x) = \prod_{j=1}^n \e^{H_j(\overline G_{j+1,n}(x)) - 1}
= \exp\left\{ \sum_{j=1}^n (H_j(\overline G_{j+1,n}(x)) - 1 ) \right\},
$$
which is easier to handle, because of its exponential form. This is
a kind of accompanying law of $X_n$. So for proving a limit theorem
we have to check the following two conditions:
\begin{itemize}
\item[(a)] $F_n (x) - \tilde F_n(x) \to 0$ for all $x \in (0,1)$, and
\item[(b)] the limit $\lim_{n \to \infty} \tilde F_n(x)$ exists for
all $ x \in (0,1)$, and the limit function is a generating function.
\end{itemize}

Also note the intuitively clear fact, that
\begin{equation} \label{rho}
\overline G^\prime_{j+1,n}(1) = \rho_{j+1} \rho_{j+2} \cdots \rho_n = \rho_{[j,n]}.
\end{equation}
This is because $\overline G^\prime_{n,n}(1) = G_n^\prime(1) = \rho_n$ and
$$
\overline G^\prime_{j,n}(1) = G_j^\prime (1) \overline G^\prime_{j+1,n}(1)
= \rho_j \overline G^\prime_{j+1,n}(1),
$$
and so by induction we have \eqref{rho}.

\subsection{Poisson and compound Poisson limits}

In this subsection we investigate the case, when the second derivative of
the offspring generating function  goes to 0 so fast that it does not appear
in the limit distribution. So the results here are the analogs of the ones
in the case of Bernoulli offsprings in \cite{GyIPV}.

In the followings Poisson$(\lambda)$ stands for a Poisson distribution
with parameter $\lambda$ if $\lambda > 0$, and Poisson$(0)$ is the degenerate distribution
at 0, while Bernoulli$(p)$ is a Bernoulli distribution with parameter $p > 0$.

The first theorem deals with the case when the immigration has
Bernoulli distribution.

\begin{theorem} \label{th1-gyipv}
Let $\{ X_n \}_{n \in \N}$ be a Galton--Watson process with immigration,
with $\varepsilon_n \sim \ \mathrm{Bernoulli}(m_{n,1})$. Assume that
\begin{itemize}
\item[(i)] $\rho_n < 1$, $\lim_{n \to \infty} \rho_n = 1$,
$\sum_{n=1}^\infty ( 1 - \rho_n) = \infty$ and
$\lim_{n \to \infty} \frac{G_n^{\prime \prime}(1)}{1 - \rho_n}= 0$,
\item[(ii)] $\lim_{n \to \infty} \frac{m_{n,1}}{1 - \rho_n} = \lambda$.
\end{itemize}
Then
$$
X_n \stackrel{\cal D}{\longrightarrow} \mathrm{Poisson}(\lambda).
$$
\end{theorem}
\medskip

Next we turn to general immigration distributions.
Let us denote the factorial  moments of the immigration distribution by
\begin{equation} \label{m_n,k}
m_{n,k} := \E [ \varepsilon_n (\varepsilon_n - 1) \ldots (\varepsilon_n - k + 1) ]
= H_n^{(k)}(1).
\end{equation}
(Clearly, this is consistent with the notation Bernoulli$(m_{n,1})$.)

The analog of Theorem 2 \cite{GyIPV} can be shown: the process has Poisson limit
even if the immigration distribution is not Bernoulli, only `close' to it.
Since the proof is also the same as there, we omit it.

\begin{theorem} \label{th2-gyipv}
Let $\{ X_n \}_{n \in \N}$ be a Galton--Watson process with immigration. Assume that
\begin{itemize}
\item[(i)] $\rho_n < 1$, $\lim_{n \to \infty} \rho_n = 1$,
$\sum_{n=1}^\infty ( 1 - \rho_n) = \infty$ and
$\lim_{n \to \infty} \frac{G_n^{\prime \prime}(1)}{1 - \rho_n}= 0$,
\item[(ii)] $\lim_{n \to \infty} \frac{m_{n,1}}{1 - \rho_n} = \lambda$,
$\lim_{n \to \infty} \frac{m_{n,2}}{1 - \rho_n}=0$.
\end{itemize}
Then
$$
X_n \stackrel{\cal D}{\longrightarrow} \mathrm{Poisson}(\lambda).
$$
\end{theorem}
\medskip

Before the following results we recall that for a finite measure $\mu$ on
$\mathbb{Z}_+ = \{ 1, 2, \ldots \}$ the \textit{compound Poisson distribution}
CP$(\mu)$ with intensity measure $\mu$ is the distribution which has
generating function
$$
\exp \left\{  \sum_{j=1}^\infty \mu\{ j \} ( x^j - 1 ) \right \},
\quad x \in [0,1].
$$

The analogs of Theorem 4 and 5 in \cite{GyIPV} are also true. 
The proof of Theorem \ref{th4-gyipv} is the same as of Theorem 4 in \cite{GyIPV}, and
also follows from the stronger Theorem \ref{th5-gyipv}, so we skip it.

\begin{theorem}  \label{th4-gyipv}
Let $\{ X_n \}_{n \in \N}$ be a Galton--Watson process with immigration. Assume that
\begin{itemize}
\item[(i)] $\rho_n < 1$, $\lim_{n \to \infty} \rho_n = 1$,
$\sum_{n=1}^\infty ( 1 - \rho_n) = \infty$ and
$\lim_{n \to \infty} \frac{G_n^{\prime \prime}(1)}{1 - \rho_n}= 0$,
\item[(ii)] $\lim_{n \to \infty} \frac{m_{n,j}}{j(1 - \rho_n)} = \lambda_j$ for
$j=1,2, \ldots, J$ with $\lambda_J=0$.
\end{itemize}
Then
$$
X_n \stackrel{\cal D}{\longrightarrow} \mathrm{CP}(\mu),
$$
where $\mu$ is a finite measure on $\{1,2, \ldots, J-1\}$ given by
$$
\mu\{ j \} = \frac{1}{j!} \sum_{i=0}^{J - j - 1} \frac{(-1)^i}{i!} \lambda_{j+i},
\quad j=1,2, \ldots, J-1.
$$
\end{theorem}
\medskip

The following theorem is more general than Theorem 5 in \cite{GyIPV} even
in the case of Bernoulli offspring distributions. And since
the proof of Theorem 5 in \cite{GyIPV} is only given by Poisson approximation we give
the analytical proof.

\begin{theorem}  \label{th5-gyipv}
Let $\{ X_n \}_{n \in \N}$ be a Galton--Watson process with immigration. Assume that
\begin{itemize}
\item[(i)] $\rho_n < 1$, $\lim_{n \to \infty} \rho_n = 1$,
$\sum_{n=1}^\infty ( 1 - \rho_n) = \infty$ and
$\lim_{n \to \infty} \frac{G_n^{\prime \prime}(1)}{1 - \rho_n}= 0$,
\item[(ii)] $\lim_{n \to \infty} \frac{m_{n,j}}{j(1 - \rho_n)} = \lambda_j$ for
$j=1,2, \ldots$, and $\limsup_{n \to \infty} \sqrt[n]{\lambda_n/n!} \leq 1$.
\end{itemize}
Then
$$
X_n \stackrel{\cal D}{\longrightarrow} Y,
$$
where the random variable $Y$ has generating function
$$
\E x^Y =
\exp \left\{ \sum_{l=1}^\infty \frac{(x-1)^l }{l!} \lambda_l \right\}.
$$
\end{theorem}
\medskip

Note that if $\lambda_n = 0$ for some $n \geq 2$, then $\lambda_m = 0$ for all $m \geq n$.
This simple fact follows from the previous two theorems.

Also notice that the assumption $\limsup_{n \to \infty} \sqrt[n]{\lambda_n/n!} \leq 1$ implies
that the limiting generating function exists for $ x \in (0,1)$. This assumption is weaker than
the one in \cite{GyIPV}, which is that for all $j = 1,2, \ldots $ the limit
\begin{equation} \label{mu}
\mu\{ j \} = \frac{1}{j!} \sum_{i=0}^{\infty} \frac{(-1)^i}{i!} \lambda_{j+i}
\end{equation}
exists. However, under this assumption the limit turns out to be a compound Poisson
random variable with intensity measure $\mu$, that is
$$
\exp \left\{ \sum_{l=1}^\infty \frac{(x-1)^l }{l!} \lambda_l \right\}
= \exp \left\{ \sum_{l=1}^\infty \mu\{ l \} (x^l - 1) \right\}. 
$$
This follows easily by Abel's theorem.

The limit can be compound Poisson even if condition \eqref{mu} fails.
If $\lambda_n = (n -1 )!/n$,
then condition \eqref{mu} does not hold even for $j = 2$, but
$$
\sum_{l=1}^\infty \frac{(x-1)^l }{l^2} =
- \int_0^{x-1} \frac{ \log ( 1- u)}{u} \d u =
\int_1^x \frac{\log ( 2- y)}{1-y} \d y
= \sum_{j=1}^\infty \mu\{ j \} (x^j - 1),
$$
with
$$
\mu \{ j \} = \frac{1}{j} \left[ \log 2 - \sum_{k=1}^{j-1} \frac{1}{k 2^k} \right]
$$
that is the limit is compound Poisson with intensity measure $\mu$.
To see that the sequence $\{ \lambda_n = \frac{(n-1)!}{n} \}_{n=1}^\infty$
can be a limit in condition (ii) consider the following example:
$\rho_n = 1- n^{-1}$ and $H_n(x) = 1 + (H(x) - 1)/n$ with the generating
function $H(x) = 1 - \log ( 2 - x)$.

We do not know whether the limit in the previous theorem necessarily has
a compound Poisson distribution.

\subsection{Negative binomial limits}

In the followings we investigate the case, when
$G_n^{\prime \prime}(1) / ( 1 - \rho_n) \not\to 0$.
In contrast to the previous subsection we show that in this case the
second derivatives do appear in the limit. Here we restrict ourselves
to the case when the immigration distribution is close to a Bernoulli distribution.

Note that for condition (a)  it is only needed
that $\rho_n < 1$, $\rho_n \to 1$ and $\sum_{n=1}^\infty ( 1 - \rho_n) = \infty$
(see the proof of Theorem 1).
Since we always assume these conditions we concentrate on condition (b),
i.e.~--in case of Bernoulli immigration-- on the existence of the limit
\begin{equation} \label{sumG}
\lim_{n \to \infty} \sum_{j=1}^n m_{j,1} ( \overline G_{j+1,n}(x) - 1).
\end{equation}

We try to compute the derivatives at 1 of the components in the sum above.
First note that since $\overline G_{j+1,n}^\prime (1) = \rho_{[j,n]}$, the
first derivative in \eqref{sumG} is
$$
\sum_{j=1}^n m_{j,1} \rho_{[j,n]} =
\sum_{j=1}^n \frac{m_{j,1}}{ 1- \rho_j} (1 - \rho_j) \rho_{[j,n]},
$$
and since $\{ (1 - \rho_j) \rho_{[j,n]} \}$ form a Toeplitz matrix, we obtain that
in order to get a limit the asymptotic order of $m_{n,1}$ must be $1 - \rho_n$.
For the second derivative we have the recursion
\begin{eqnarray}  \label{der2G}
\overline G_{j,n}^{\prime \prime}(x)
& = & \frac{\d^2}{\d x^2} G_j(\overline G_{j+1,n}(x)) =
\frac{\d}{\d x} \left( G_j^\prime (\overline G_{j+1,n}(x)) 
\overline G_{j+1,n}^\prime(x) \right) \nonumber \\
& = & G_j^{\prime \prime} (\overline G_{j+1,n}(x)) \overline G_{j+1,n}^\prime(x)^2 +
G_j^\prime (\overline G_{j+1,n}(x))  \overline G_{j+1,n}^{\prime \prime}(x).
\end{eqnarray}
Substituting $x= 1$ for $j \leq n$
(recall that $\overline G_{n+1,n}(x) = x$) by \eqref{rho} we obtain that
\begin{equation} \label{der2G(1)}
\overline G_{j,n}^{\prime \prime}(1) = G_j^{\prime \prime}(1) \rho_{[j,n]}^2 +
\rho_j  \overline G_{j+1,n}^{\prime \prime}(1),
\end{equation}
and $\overline G_{n,n}^{\prime \prime}(1)= G_{n}^{\prime \prime}(1)$.
Induction argument shows that
$$
\overline G_{j+1,n}^{\prime \prime} (1) = \sum_{i=j+1}^n G_i^{\prime \prime}(1)
\rho_{[j,i-1]} \rho_{[i,n]}^2,
$$
so the second derivative in \eqref{sumG} is
$$
\sum_{j=1}^n m_{j,1} \overline G_{j+1,n}^{\prime \prime} (1)
= \sum_{i=2}^n G_i^{\prime \prime} (1) \rho_{[i,n]}^2
\sum_{j=1}^{i-1} m_{j,1} \rho_{[j, i- 1]}.
$$
As we have seen $m_n \approx 1 - \rho_n$, and so to get a limit for the second derivative
we must have $G_n^{\prime \prime}(1) \approx 1 - \rho_n$.
These heuristic argument kind of shows the necessity of the assumptions in the following
theorem. We also note that Step 2 of the proof Theorem \ref{th-nb} shows that
the higher derivatives must have the same order, namely $1 - \rho_n$. We are not
able to calculate these contributions, therefore we have to assume condition (iii). 

The reasoning above can be made rigorous to show the following: If
$\lim_{n \to \infty} G_n^{\prime \prime}(1) = a$ exists and the inhomogeneous Galton--Watson
process with immigration has a \emph{proper} limit distribution
with finite second moment, then
necessarily $a = 0$. This immediately implies that $\lim_{n \to \infty} G_n (x) = x$,
i.e.~there is no critical nontrivial branching mechanism with finite second moment,
which can cause a proper limit distribution with finite second moment.
This result is in complete accordance
with Theorem 1 (ii) by Foster and Williamson \cite{foster} in the homogeneous case.

We will show that in this setup the limit distribution is the negative binomial distribution.
A random variable $X$ has negative binomial distribution with parameters $r >0$ and 
$p \in (0,1)$, denoted by NB$(r,p)$, if $\p \{ X = k \} = \binom{k+r - 1}{r-1} (1 - p)^r p^k$,
$k=0,1,2, \ldots$, where the binomial coefficient is defined by
$\binom{k+r - 1}{r-1} = \frac{( k+r-1) (k+r -2 ) \cdots r}{k!}$. The generating function is
$$
\E x^X = \left( \frac{1 - p}{1 - p x} \right)^r.
$$

The following theorem holds.

\begin{theorem} \label{th-nb}
Let $\{ X_n \}$ be a Galton--Watson process with immigration, with general offspring and
immigration distribution, such that the followings hold:
\begin{itemize}
\item[(i)] $\rho_n < 1$, $\lim_{n \to \infty} \rho_n = 1$,
$\sum_{n=1}^\infty ( 1 - \rho_n) = \infty$,

\item[(ii)] $\lim_{n \to \infty} \frac{G_n^{\prime \prime}(1)}{ 1 - \rho_n} =
\nu \in (0, \infty)$,

\item[(iii)] $\lim_{n \to \infty}
\frac{G_n^{(s)}(1)}{1 - \rho_n} =0$, for all $s \geq 3$,

\item[(iv)] $\lim_{n \to \infty} \frac{m_{n,1}}{1 - \rho_n} = \lambda$ and
$\lim_{n \to \infty} \frac{m_{n,2}}{1 - \rho_n} = 0$.
\end{itemize}
Then
$$
X_n \stackrel{\cal D}{\longrightarrow} \mathrm{NB}
(2 \lambda / \nu , \nu/ (2 + \nu)).
$$
\end{theorem}
\medskip

As we already mentioned, assumption (ii) means that the second derivatives have
a significant role in the limit, while assumption (iii) ensures that the third
and higher derivatives do not count. The proof of the theorem basically lays on
determining the asymptotic of $\overline G_{j+1,n}^{(k)}$, combined with a
relatively closed formula for the coefficient of $f^{\prime \prime }(g)$ in
$\frac{\d^k}{\d x^k} f(g(x))$, see \eqref{comp-der}. This `shows' that
to handle the case, when the higher derivatives also count, a different approach is
needed, or at least the calculations become more technical. It would be also
interesting to extend the results to more general immigration distribution, that is
to know whether a kind of analog of Theorem \ref{th5-gyipv} remains true.

Finally, we note that if $\lambda > 0$ is fixed and $\nu \to 0$, then
$$
\lim_{ \nu \to 0} \mathrm{NB}(2 \lambda / \nu,  \nu / ( 2 + \nu) )
\to \mathrm{Poisson}( \lambda).
$$
It is easy to check that the proof of the theorem remains correct in this case,
and we obtain that Theorem \ref{th-nb} holds if $\nu = 0$ and in this case the
limit distribution is Poisson$( \lambda )$, which is exactly the statement
of Theorem \ref{th2-gyipv}.

\subsection{The case $\sum_{n=1}^\infty ( 1 - \rho_n ) < \infty$}

Let us consider the nearly critical inhomogeneous Galton--Watson process with general
offspring distribution
and general immigration. It was always assumed that $\rho_n$ does not
converge too fast to 1, that is $\sum_{n=1}^\infty ( 1 - \rho_n) = \infty$.
In this section we
investigate the case, when this assumption does not hold, and assume that
$$
\prod_{n=2}^\infty \rho_n = \rho \in (0,1).
$$
(Note that this is equivalent to the assumption $\sum_{n=1}^\infty ( 1 -\rho_n) < \infty$.)
We show that in this case the limit distribution exists under very general assumptions, but
the process does not contain enough randomness as in the previous cases: the limit
distribution explicitly contains each offspring- and immigration distribution.
That is $\rho_n$ tends to 1 so fast, that the process has no time
to forget the initial distribution. This shows that the right assumption is indeed
$ \sum_{n=1}^\infty ( 1 - \rho_n) = \infty$, as it was investigated in \cite{GyIPV}.
In fact, the case treated in this section is much simpler.

We state  the theorem in the most general setup.

\begin{theorem} \label{th-gyors-konv}
Let $X_n$ be the  Galton--Watson process with immigration, described above.
If $\prod_{n=2}^\infty \rho_n = \rho \in (0,1)$ and $\sum_{n=1}^\infty m_{n,1} < \infty$,
then
$$
X_n \stackrel{{\cal D}}{\longrightarrow} Y,
$$
where $Y$ has generating function
$$
g(x) = \prod_{j=1}^\infty  H_j( \overline G_{j+1, \infty}(x)),
$$
with
$$
\overline G_{j+1,\infty}(x) = \lim_{n \to \infty }  \overline G_{j+1, n}(x).
$$
\end{theorem}

Let us see an example in the simplest case, when both the offspring distribution and the
immigration distribution is Bernoulli. In this case we can compute the limit generating
function.

\begin{example}
Let $\rho_n = 1 - \frac{1}{n^2}$. Clearly
$\sum_{n=1}^\infty ( 1 - \rho_n) = \pi^2/6 < \infty$. For the product we have
\begin{eqnarray*}
\rho_{[1,n]}
& = & \rho_2 \rho_3 \ldots \rho_n =
\left(1 - \frac{1}{2^2} \right) \left(1 - \frac{1}{3^2} \right)
\cdots \left( 1 - \frac{1}{n^2} \right) \\
& = & \frac{(2^2 - 1) (3^2 - 1) \ldots (n^2 - 1)}{(n!)^2} \\
& = & \frac{1 \cdot 3 \cdot 2 \cdot 4 \cdot \ldots \cdot (n-1) \cdot (n+1) }{(n!)^2} \\
& = & \frac{n+1}{2 n} \to \frac{1}{2} = \rho.
\end{eqnarray*}
In this case $\overline G_{j+1,n}(x) = 1 + \rho_{[j,n]}(x - 1)$, so
$\overline G_{j+1, \infty}(x) = 1 + \rho (x - 1) / \rho_{[1,j]}.$
Since $H_j(x)= 1 - m_{j,1}( 1 - x)$, for the limit generating function we have
$$
g(x) = \prod_{j=1}^\infty
\left[ 1 - \frac{m_{j,1}}{\rho_{[1,j]}} \frac{1}{2} ( 1 - x ) \right] =
\prod_{j=1}^\infty
\left[ 1 - m_{j,1} \frac{j}{j+1} ( 1 - x ) \right].
$$
Choose $m_{j,1} = (j+1) / j^3$, so
$$
g(x) = \prod_{j=1}^\infty
\left[ 1 - \frac{ 1 - x}{j^2}  \right] = \frac{\sin ( \pi \sqrt{1- x} )}{\pi \sqrt{1- x}}.
$$
\end{example}

The following example is Example 1 in \cite{GyIPV}, which shows that in a special case
the limit can be Poisson even in this setup.

\begin{example}
Assume that the offspring distribution is Bernoulli$(\rho_n)$
and the immigration $\varepsilon_n$ has Poisson$(m_{n,1})$ distribution, where
$ \sum_{n=1}^\infty ( 1 - \rho_n) < \infty$ and $\sum_{n=1}^\infty m_{n,1} < \infty$.
As before $\overline G_{j+1, \infty}(x) = 1 + \rho (x - 1) / \rho_{[1,j]},$ and
since $H_n(x) = \e^{m_{n,1} (x - 1)}$ we have that
$$
g(x) = \exp \left\{ (x - 1) \sum_{n=1}^\infty \frac{m_{n,1} \rho}{\rho_{[1,n]}} \right\},
$$
and since the sum in the exponent is finite, this is Poisson distribution with mean
$\sum_{n=1}^\infty \frac{m_{n,1} \rho}{\rho_{[1,n]}}$.
\end{example}

\subsection{The linear fractional case}

The importance of the linear fractional generating functions in branching
processes is that, that this is basically the only example when explicit
computation can be done. In this subsection we investigate this example in
detail, which actually helped to find the general form of Theorem \ref{th-nb}.

The linear fractional generating function has the form 
$f(s) = 1 - \frac{\alpha}{1- \beta} + \frac{\alpha s}{1 - \beta s},$
where $\alpha, \beta \in (0,1)$, $\alpha + \beta \leq 1$.
For the first two derivatives we have
\begin{eqnarray*}
f^\prime(1) & = & \frac{\alpha}{ (1-\beta)^2} \\
f^{\prime \prime}(1) & = & \frac{2 \alpha \beta}{(1 - \beta)^3},
\end{eqnarray*}
and these determine the parameters $\alpha, \beta$:
\begin{eqnarray*}
\alpha &  = & \frac{4 f^\prime(1)^3}{(2 f^\prime(1) + f^{\prime \prime}(1))^2} \\
\beta & = & \frac{f^{\prime \prime}(1)}{2 f^\prime(1) + f^{\prime \prime}(1)}.
\end{eqnarray*}

If $f_1(s)$ and $f_2(s)$ are both linear fractional generating
functions then so is $f_1(f_2(s))$, that is if $G_n$ is linear fractional for all $n$,
so is $\overline G_{j,n}(s)$ for all $j,n$, with parameters
\begin{eqnarray*}
\alpha_{j,n} & = & \frac{ 4 \overline G_{j,n}^{\prime}(1)^3}
{( 2 \overline G_{j,n}^{\prime}(1) + \overline G_{j,n}^{\prime \prime}(1))^2}, \\
\beta_{j,n} & = & \frac{\overline G_{j,n}^{\prime \prime }(1)}
{ 2 \overline G_{j,n}^{\prime}(1) + \overline G_{j,n}^{\prime \prime}(1)}.
\end{eqnarray*}
As we have seen $\overline G_{j,n}^\prime(1) = \rho_{[j-1,n]}$ and
$\overline G_{j,n}^{\prime \prime}(1) = \sum_{i=j}^n G_i^{\prime \prime}(1)
\rho_{[j-1,i-1]} \rho_{[i,n]}^2$.
The generating function has the form
$F_n(x) = \prod_{j=1}^n H_j( \overline G_{j+1,n}(x))$.
Assuming Bernoulli immigration, i.e.~$H_n(x) = 1 + m_n (x - 1)$,
the corresponding accompanying generating function is
\begin{eqnarray*}
\tilde F_n (x)
& = & \exp \left\{
\sum_{j=1}^n m_j (\overline G_{j+1,n}(x) - 1) \right\} \\
& = & \exp \left\{ \sum_{j=1}^n m_j
\left( - \frac{\alpha_{j+1,n}}{1 - \beta_{j+1,n}} +
\frac{ \alpha_{j+1,n} x}{1 - \beta_{j+1,n} x} \right)
\right\}. 
\end{eqnarray*}
Some calculation shows
$$
\frac{\alpha_{j+1,n}}{1 - \beta_{j+1,n} x}
= \frac{\rho_{[j,n]}}
{\left( 1 + \sum_{i=j+1}^n \frac{G_i^{\prime \prime}(1) \rho_{[i,n]}}{2 \rho_i} \right)
\left( 1 + \sum_{i=j+1}^n \frac{G_i^{\prime \prime}(1) \rho_{[i,n]}}{2 \rho_i}
(1 - x) \right)}.
$$

Eventually, we obtained that in the linear fractional case the limit exists
(under the assumption $\sum_{n=1}^\infty ( 1  - \rho_n) = \infty$) if and only if
$$
\lim_{n\to  \infty}
\sum_{j=1}^n  m_j \rho_{[j,n]}
\left[
\left( 1 + \sum_{i=j+1}^n \frac{G_i^{\prime \prime}(1) \rho_{[i,n]}}{2 \rho_i} \right)
\left( 1 + \sum_{i=j+1}^n \frac{G_i^{\prime \prime}(1) \rho_{[i,n]}}{2 \rho_i}
(1 - x) \right)
\right]^{-1}
$$
exists. It is easy to check that if (ii) in Theorem \ref{th-nb} holds, then so is
(iii), that is the limit is negative binomial distribution.

\section{Proofs}

We use the continuity theorem for generating functions (\cite{feller} p.~280) and
Lemmas 5 and 6 in \cite{GyIPV} without any further reference. We also
frequently use the simple facts that
\begin{equation} \label{prod-inequ}
\left| \prod_{k=1}^n z_k - \prod_{k=1}^n w_k \right| \leq
\sum_{k=1}^n | z_k - w_k |,
\end{equation}
for $z_k, w_k \in [-1,1]$, $k=1,2, \ldots, n$, (see for example
Lemma 3 in \cite{GyIPV}) and that
$|\e^u - 1 - u | \leq u^2$ for $|u| \leq 1/2$.

\begin{proof}{Theorem \ref{th1-gyipv}}
By the assumptions the immigration $\varepsilon_n$ has Bernoulli$(m_{n,1})$ distribution, that
is $H_n(x) = 1 + m_{n,1} (x - 1)$, in which case formula \eqref{gen-fgv} reduces to
$$
F_n(x) = \prod_{j=1}^n \left[ 1 + m_{j,1} ( \overline G_{j+1,n}(x) - 1 ) \right].
$$
Let us define the generating function
$$
\tilde F_n (x) = \prod_{j=1}^n \e^{m_{j,1} ( \overline G_{j+1,n}(x) - 1)}.
$$
Estimation \eqref{lb-gen-fgv} gives that if $m_{n,1} \to 0$ then
$\max_{1 \leq j \leq n} m_{j,1} |\overline G_{j+1,n}(x) - 1| \to 0$.

We want to show that for all $x \in (0,1)$
\begin{itemize}
\item[(a)] $\tilde F_n(x) - F_n(x) \to 0$, as $n \to \infty$, and
\item[(b)] $\lim_{n \to \infty} \tilde F_n(x) = \e^{\lambda (x - 1)}$.
\end{itemize}

Condition (a) is the easier to handle. Using
\eqref{prod-inequ}, the inequality $|\e^u - 1 - u | \leq u^2$ and
the mean value theorem combined with the monotonicity of the derivative, we obtain
\begin{eqnarray*}
| \tilde F_n(x) - F_n(x) |
& = &  \left| \prod_{j=1}^n \e^{m_{j,1} ( \overline G_{j+1,n}(x) - 1)} -
\prod_{j=1}^n \left[ 1 + m_{j,1} ( \overline G_{j+1,n}(x) - 1 ) \right] \right| \\
& \leq & \sum_{j=1}^n \left| \e^{m_{j,1} ( \overline G_{j+1,n}(x) - 1)}
-  \left[ 1 + m_{j,1} ( \overline G_{j+1,n}(x) - 1 ) \right] \right| \\
& \leq & \sum_{j=1}^n m_{j,1}^2 | \overline G_{j+1,n}(x) - 1 |^2 \\
& \leq & \sum_{j=1}^n m_{j,1}^2 | x - 1 |^2  \left( \overline G^\prime_{j+1,n}(1) \right)^2.
\end{eqnarray*}
Using \eqref{rho}, finally we have
$$
| \tilde F_n(x) - F_n(x) | \leq
| x - 1 |^2 \sum_{j=1}^n m_{j,1}^2  \rho_{[j,n]}^2,
$$
and by Lemma 5 in \cite{GyIPV} this goes to 0 under the assumptions of the theorem.

Now consider condition (b). This convergence is equivalent to
$$
\sum_{j=1}^n m_{j,1} [ \overline G_{j+1,n}(x) - 1 ] \to \lambda (x - 1).
$$
Let $x \in (0,1)$. The convexity of the generating
functions implies $G_n(x) \geq 1 - \rho_n + \rho_n x = 1 + \rho_n (x - 1)$. And so
$$
\overline G_{n-1,n}(x) = G_{n-1} (G_n(x)) \geq
1 + \rho_{n-1} (G_n(x) - 1) \geq 1 + \rho_{n-1} \rho_n (x - 1).
$$
Induction gives that
\begin{equation} \label{lb-gen-fgv}
\overline G_{j+ 1, n}(x) \geq 1 + \rho_{j + 1} \rho_{j+ 2} \cdots \rho_n (x - 1)
= 1 + \rho_{[j,n]} (x - 1),
\end{equation}
which implies the lower bound
$$
\sum_{j=1}^n m_{j,1} [ \overline G_{j+1,n}(x) - 1 ] \geq
\sum_{j=1}^n m_{j,1} \rho_{[j,n]} (x - 1) \to \lambda (x - 1).
$$

For the upper bound note that by the previous estimations
$\overline G_{j+1,n}(x) \in (1 - \rho_{[j,n]}, 1) $ for any $x \in (0, 1)$. Again by
convexity, for $ y \in (1 - \rho_{[j,n]}, 1)$ we have
\begin{eqnarray*}
G_j(y)
& = & G_j \left( \frac{1-y}{\rho_{[j,n]}} (1 - \rho_{[j,n]}) + 1
- \frac{1-y}{\rho_{[j,n]}} \right) \\
& \leq & \frac{1-y}{\rho_{[j,n]}} G_j(1- \rho_{[j,n]})
+  1 - \frac{1-y}{\rho_{[j,n]}} \\
& = & 1 + \vartheta_{j,n} ( y - 1),
\end{eqnarray*}
where
$$
\vartheta_{j,n} = \frac{1 - G_j( 1 - \rho_{[j,n]})}{\rho_{[j,n]}}
= G_j^\prime (\xi_{j,n} ),
$$
by the mean value theorem.
Since $\rho_{[j,n]} \to 0$ thus $\xi_{j,n} \to 1$
as $ n \to \infty$, $\vartheta_{j,n} \uparrow \rho_j$.
So we get
$$
\overline G_{j,n}(x) = G_j(\overline G_{j+1, n}(x)) \leq
1 + \vartheta_{j,n} ( \overline G_{j+1, n}(x) - 1),
$$
and so induction gives
\begin{equation} \label{ub-gen-fgv}
\overline G_{j+1,n}(x) \leq 1 + \vartheta_{j+1,n} \vartheta_{j+2,n} \cdots \vartheta_{n,n} (x - 1)
=: 1 + \vartheta_{[j,n]} (x - 1).
\end{equation}
Summarizing we have
$$
\sum_{j=1}^n m_{j,1} [ \overline G_{j+1,n}(x) - 1 ] \leq
\sum_{j=1}^n m_{j,1} \vartheta_{[j,n]} (x - 1).
$$
Note that in the case of Bernoulli offspring distributions the upper and lower bounds are equal.

So we have to check that under what conditions
$$
\sum_{j=1}^n m_{j,1} \vartheta_{[j,n]} \to \lambda.
$$
For this, exactly the same way as in Lemma 5 in \cite{GyIPV} we only need that
the sequence $b_{n,j} : = ( 1 - \rho_j) \vartheta_{[j,n]}$ form a Toeplitz matrix.
The only nontrivial condition is
$$
\sum_{j=1}^n b_{n,j} \to 1.
$$
We may write
\begin{eqnarray*}
\sum_{j=1}^n b_{n,j}
& = & ( 1 - \rho_1) \prod_{l = 2}^n \vartheta_{l,n} +
(1 - \rho_2) \prod_{l = 3}^n \vartheta_{l,n} + \cdots \\
&& \hspace{81pt} + ( 1 - \rho_{n-1}) \vartheta_{n,n} + ( 1 - \rho_n) \\
& = & 1 - \bigg[ ( \rho_n - \vartheta_{n,n}) +
(\rho_{n-1} - \vartheta_{n-1,n}) \vartheta_{n,n} + \cdots \\
&& \hspace{86pt} + (\rho_2 - \vartheta_{2,n}) \prod_{l=3}^n \vartheta_{l,n}
+ \rho_1 \prod_{l=2}^n \vartheta_{l,n}
\bigg].
\end{eqnarray*}
Note that in the bracket every term is nonnegative, since 
$\vartheta_{j,n} = G_j^\prime (\xi_{j,n}) \leq G_j^\prime (1) = \rho_j$,
because $G_j^\prime$ is monotone increasing since $G_j$ is convex.
So we need that
$$
( \rho_n - \vartheta_{n,n}) +
(\rho_{n-1} - \vartheta_{n-1,n}) \vartheta_{n,n} + \cdots +
(\rho_2 - \vartheta_{2,n}) \prod_{l=3}^n \vartheta_{l,n}
+ \rho_1 \prod_{l=2}^n \vartheta_{l,n} \to 0.
$$
By definition we have
$\vartheta_{j,n} = G_j^\prime ( \xi_{j,n} )$,
where $\xi_{j,n} \in ( 1- \rho_{[j,n]}, 1)$, and since
$\rho_j = G_j^\prime(1)$ the mean value theorem again gives
$$
\rho_j - \vartheta_{j,n} = G_j^\prime (1) - G_j^\prime ( \xi_{j,n})
= (1 - \xi_{j,n}) G_j^{\prime \prime}(\xi_{j,n}^\prime) \leq
\rho_{[j,n]} G_j^{\prime \prime}(1),
$$
where $\xi_{j,n}^\prime \in (\xi_{j,n}, 1)$. Therefore we obtain
$$
\sum_{j=1}^n ( \rho_j - \vartheta_{j,n}) \vartheta_{[j,n]}
\leq \sum_{j=1}^n \rho_{[j,n]} G_j^{\prime \prime}(1) \vartheta_{[j,n]}
\leq \sum_{j=1}^n \frac{G_j^{\prime \prime}(1)}{1 - \rho_j}
( 1 - \rho_j) \rho^2_{[j,n]},
$$
which, according to Lemma 5 \cite{GyIPV}, goes to 0, if
$\frac{G_n^{\prime \prime}(1)}{1 - \rho_n} \to 0$.
The proof is  ready.
\end{proof}

\begin{proof}{Theorem \ref{th5-gyipv}}
By formula \eqref{gen-fgv} for the generating function we have to show that
for all $x \in (0,1)$
$$
F_n(x) = \prod_{j=1}^n H_j(\overline G_{j+1,n}(x)) \to 
\exp \left\{ \sum_{l=1}^\infty \frac{(x-1)^l }{l!} \lambda_l \right\}.
$$
Let us define the function
$$
\tilde F_n(x) = \prod_{j=1}^n \e^{H_j(\overline G_{j+1,n}(x)) - 1}.
$$
Using the estimation \eqref{lb-gen-fgv} and that $0 \leq 1 - H_j(x) \leq m_{j,1} (1- x)$,
we have
\begin{eqnarray*}
| F_n(x) - \tilde F_n (x) |
& \leq & \sum_{j=1}^n  \left| \e^{H_j(\overline G_{j+1,n}(x)) - 1} -
H_j(\overline G_{j+1,n}(x))  \right| \\
& \leq & \sum_{j=1}^n \left( H_j(\overline G_{j+1,n}(x)) - 1 \right)^2 \\
& \leq & \sum_{j=1}^n m_{j,1}^2 \rho_{[j,n]}^2 \to 0, 
\end{eqnarray*}
since $m_{j,1} / ( 1 - \rho_j) \to \lambda_1$ implies $m_{j,1}^2 / ( 1 - \rho_j) \to 0$.
So we obtain the convergence $F_n(x) - \tilde F_n(x) \to 0$, for all $x \in (0,1)$.
Therefore what we have to show is that
$$
\sum_{j=1}^n \left[ H_j(\overline G_{j+1,n}(x)) - 1 \right] \to 
\sum_{l=1}^\infty \frac{(x-1)^l}{l!} \lambda_l.
$$

Let us fix an $\varepsilon > 0$ and an $x \in (0,1)$. There is an $l_0$ such that
\begin{equation} \label{l_0}
\left| \sum_{l=l_0 + 1}^\infty \frac{(x-1)^l}{l!} \lambda_l \right| \leq \varepsilon \quad
\textrm{and } \ 
\frac{(1 -x)^{l_0 + 1}}{(l_0 +1)!} \lambda_{l_0 + 1} \leq \varepsilon.
\end{equation}
By Lemma 6 in \cite{GyIPV} we have
$$
H_j(x) = \sum_{l=0}^{l_0} \frac{m_{j,l}}{l!} (x - 1)^l + R_{j,l_0 + 1}(x),
$$
where $|R_{j,l_0 + 1}(x)| \leq (1 - x)^{l_0 + 1} m_{j,l_0 +  1} / (l_0 + 1)!$. Thus
\begin{eqnarray*}
&& \sum_{j=1}^n \left[ H_j(\overline G_{j+1,n}(x)) - 1 \right] \\
&& =  \sum_{j=1}^n \left\{ \sum_{l=1}^{l_0} \frac{m_{j,l}}{l!} (\overline G_{j+1,n}(x) - 1)^l
+ R_{j,l_0 + 1}(\overline G_{j+1,n}(x)) \right\} \\
&& =  \sum_{l=1}^{l_0} \frac{(x - 1)^l}{l!} \sum_{j=1}^n m_{j,l} \rho_{[j,n]}^l
+ \sum_{j=1}^n R_{j,l_0+1} ( \overline G_{j+1,n}(x) ) \\
&& \phantom{=} + \sum_{l=1}^{l_0} \sum_{j=1}^n \frac{m_{j,l}}{l!}
\left[ \left( \overline G_{j+1,n}(x) - 1 \right)^l - \rho_{[j,n]}^l (x - 1)^l \right] \\
&& =    \sum_{l=1}^{l_0} \frac{(x - 1)^l}{l!} \sum_{j=1}^n m_{j,l} \rho_{[j,n]}^l +
I_1 + I_2.
\end{eqnarray*}
Using \eqref{lb-gen-fgv} again
\begin{eqnarray*}
| I_1 |
& \leq & \sum_{j=1}^n \frac{m_{j,l_0 +1}}{(l_0 + 1)!}
| \overline G_{j+1,n}(x) - 1 |^{l_0 + 1} \\
& \leq & \sum_{j=1}^n \frac{m_{j,l_0 +1}}{(l_0 + 1)!} \rho_{[j,n]}^{l_0 + 1} ( 1 - x)^{l_0 + 1}
\to \frac{( 1 - x)^{l_0 + 1}}{(l_0 + 1)!} \lambda_{l_0 + 1}
\end{eqnarray*}
by Lemma 5 in \cite{GyIPV},
thus \eqref{l_0}, the choice of $l_0$, shows that $| I_1 | \leq 2 \varepsilon$ for $n$ large
enough. In order to estimate $I_2$ we use the inequalities \eqref{lb-gen-fgv} and
\eqref{ub-gen-fgv} we have
\begin{eqnarray*}
&& \left| \left( \overline G_{j+1,n}(x) - 1 \right)^l - \rho_{[j,n]}^l (x - 1)^l \right| \\
&& 
\leq \left| \left( \overline G_{j+1,n}(x) - 1 \right) - \rho_{[j,n]} (x - 1) \right|
l \rho_{[j,n]}^{l-1}  ( 1 - x)^{l-1} \\
&&  \leq   l ( 1 - x)^l ( \rho_{[j,n]} - \vartheta_{[j,n]}) \rho_{[j,n]}^{l-1},
\end{eqnarray*}
therefore
$$
| I_2 | \leq \sum_{l=1}^{l_0} \frac{(1-x)^l}{(l-1)!}
\sum_{j=1}^n m_{j,l}  ( \rho_{[j,n]} - \vartheta_{[j,n]}) \rho_{[j,n]}^{l-1},
$$
which goes to 0, due to Lemma 5 in \cite{GyIPV} and to our assumptions.

Finally for the difference
\begin{eqnarray*}
&& \left| \sum_{j=1}^n \left[ H_j(\overline G_{j+1,n}(x)) - 1 \right] - 
\sum_{l=1}^\infty \frac{(x-1)^l}{l!} \lambda_l \right| \\
&& \leq \left| \sum_{l=1}^{l_0} \frac{(x - 1)^l}{l!}
\left( \sum_{j=1}^n m_{j,l} \rho_{[j,n]}^l - \lambda_l \right) \right|
+ \left| \sum_{l=l_0 +1}^\infty \frac{(x-1)^l}{l!} \lambda_l \right|
+ | I_1 | + |I_2|,
\end{eqnarray*}
where the first and fourth term converge to 0, while the sum of second and the third is less
than $3 \, \varepsilon$. Since $\varepsilon$ was arbitrary the proof is ready.
\end{proof}

\begin{proof}{Theorem \ref{th-nb}}
We separate three cases. First we assume that the offspring generating functions are second
degree polynomials, that is each particle has at most two offsprings. Then we extend the proof
when for the higher derivatives only assumption (iii) holds. In this two cases we
assume Bernoulli immigration. In the third case we extend the result for general immigration.
Step 1 shows the main idea without the technical difficulty to handle the higher derivatives.

\noindent \textbf{Step 1.}
Let us assume that
\begin{equation} \label{G-assumption}
\textrm{deg } G_n = 2 \quad \textrm{for all } n \in  \N.
\end{equation}

First we obtain a recursion like \eqref{der2G(1)}.
What we need is a general form of a derivative of a composite function. Let
$$
h_k(x) = \frac{\d^k}{\d x^k} f(g(x)).
$$
According to Lemma 5.6 in  \cite{petrov} the general formula for the derivatives of
a composite function is
\begin{equation} \label{comp-der}
\frac{\d^k}{\d x^k} f(g(x)) =
k! \sum_{s=1}^k f^{(s)}(g(x))
\sum_{(\nu_1, \ldots, \nu_k)} \prod_{m=1}^k
\frac{1}{\nu_m !} \left( \frac{1}{m!} g^{(m)}(x) \right)^{\nu_m},
\end{equation}
where the summation is carried out over all nonnegative integer solutions
of the equation system:
\begin{eqnarray*}
&&\nu_1 + 2 \nu_2 + \cdots + k \nu_k = k \\
&&\nu_1 + \nu_2 + \cdots + \nu_k = s.
\end{eqnarray*}
Clearly the coefficient of $f^\prime(g)$ in $h_k$ is $g^{(k)}$, and
induction argument shows that
$$
\textrm{coefficient of } f^{\prime \prime}(g) \textrm{ in } h_k =
\left\{
\begin{array}{ll}
\sum_{i=1}^{\frac{k-1}{2}} \binom{k}{i} g^{(i)} \, g^{(k-i)},
&  k \textrm{ odd}, \vspace{4pt} \\
\sum_{i=1}^{\frac{k}{2}-1} \binom{k}{i} g^{(i)} \, g^{(k-i)}
+ \frac{1}{2} \binom{k}{ k/2} \big( g^{(\frac{k}{2})} \big)^2,
&  k \textrm{ even}.
\end{array}
\right.
$$
For simplicity introduce
$$
a_{k,i}= \left\{
\begin{array}{ll}
\binom{k}{i}, & \textrm{if } i < \frac{k}{2}, \\
\frac{1}{2} \binom{k}{k/2}, & \textrm{if } i = \frac{k}{2}.
\end{array}
\right.
$$
Apply this result to $\overline G_{j,n}(x) = G_j( \overline G_{j+1,n}(x))$
for $j \leq n$, and note that the  third and higher derivatives vanishes.
Substituting $x=1$ we have
\begin{equation} \label{derkG(1)}
\overline G_{j,n}^{(k)}(1) =
\rho_j \overline G_{j+1,n}^{(k)}(1) +
G_j^{\prime \prime}(1) \sum_{i=1}^{k/2} a_{k,i} \overline G_{j+1,n}^{(i)}(1)
\overline G_{j+1,n}^{(k-i)}(1).
\end{equation}
From this, easy induction argument gives the following:
\begin{equation} \label{k-th der-formula}
\overline G_{j+1,n}^{(k)}(1) =
\sum_{i=1}^{k/2} a_{k,i} \sum_{l=j+1}^n \rho_{[j,l-1]}
G_l^{\prime \prime}(1) \overline G_{l+1,n}^{(i)}(1)
\overline G_{l+1,n}^{(k-i)}(1).
\end{equation}

Until now everything hold in general. Now we use our assumptions to prove
\begin{equation} \label{k-th der-asympt}
\overline G_{j+1,n}^{(k)} ( 1 ) = \left( \frac{k!}{2^{k-1}} \nu^{k-1} + o(1) \right)
\rho_{[j,n]} ( 1 - \rho_{[j,n]} )^{k-1}
+ \rho_{[j,n]} o(1),
\end{equation}
where $o(1) \to 0$ as $j, n \to \infty$.
The proof goes by induction, uses the recursion
\eqref{k-th der-formula} and the identity
$$
\frac{(k+1)!}{2^k} = \frac{1}{k} \sum_{i=1}^{(k+1)/2}
a_{k+1,i} \frac{i!}{2^{i-1}} \frac{(k+1 - i)!}{2^{k-i}}.
$$
Formula \eqref{k-th der-asympt} is  true for $k=1,2$, and let us assume
that it is true until some $k \geq 2$.
Then using the induction hypothesis
\begin{eqnarray*}
\overline G_{j+1,n}^{(k+1)}(1)
& = & \sum_{i=1}^{\frac{k+1}{2}} a_{k+1,i} \sum_{l=j+1}^n \rho_{[j,l-1]}
G_l^{\prime \prime}(1) \overline G_{l+1,n}^{(i)}(1)
\overline G_{l+1,n}^{(k+1-i)}(1) \\
& = & \sum_{i=1}^{\frac{k+1}{2}} a_{k+1,i} \sum_{l=j+1}^n \rho_{[j,l-1]}
G_l^{\prime \prime}(1) \\
& & \times
\left[  \left( \frac{i!}{2^{i-1}} \nu^{i-1} + o(1) \right)
\rho_{[l,n]} ( 1 - \rho_{[l,n]} )^{i-1}
+ \rho_{[l,n]} o(1) \right] \\
&& \times \left[  \left( \frac{(k - i + 1)!}{2^{k-i}} \nu^{k-i} + o(1) \right)
\rho_{[l,n]} ( 1 - \rho_{[l,n]} )^{k-i}
+ \rho_{[l,n]} o(1) \right] \\
& = &  \sum_{i=1}^{\frac{k+1}{2}} a_{k+1,i}
\rho_{[j,n]} \sum_{l=j+1}^n \frac{G_l^{\prime \prime}(1)}{\rho_l ( 1 - \rho_l)}
( 1 - \rho_l ) \rho_{[l,n]} \\
& & \times
\left[ \left( \frac{(k+i-1)! i!}{2^{k-1}} \nu^{k-1} + o(1) \right)
(1 - \rho_{[l,n]})^{k-1} + o(1)
\right] \\
&  = & \left( \frac{(k+1)!}{2^{k}} \nu^{k} + o(1) \right)
\rho_{[j,n]} ( 1 - \rho_{[j,n]} )^{k}
+ \rho_{[j,n]} o(1),
\end{eqnarray*}
where we used the estimation
$$
\sum_{l= j+1}^n ( 1 - \rho_l) \rho_{[l,n]} ( 1 - \rho_{[l,n]})^k =
\frac{(1 - \rho_{[j,n]})^{k+1}}{k+1} +
O ( \max_{j \leq l \leq n} ( 1 - \rho_l ) \rho_{[l,n]}).
$$
Note that the left side is a Riemann approximation of the integral
$\int ( 1 - y)^k \d y$, corresponding to the partition
$\{ \rho_{[l,n]} \}_{l= j+1}^n$.
So \eqref{k-th der-asympt} is proved.

From \eqref{k-th der-asympt} using  also assumption (iv) in the theorem
we have
\begin{align*}
\sum_{j=1}^n m_{j,1} \overline G_{j+1,n}^{(k)}(1)
 = & 
\sum_{j=1}^n \frac{m_{j,1}}{1 - \rho_j}
\left( \frac{k!}{2^{k-1}} \nu^{k-1} + o(1) \right)
\rho_{[j,n]} ( 1 - \rho_j)   ( 1 - \rho_{[j,n]})^{k-1} \\
&   + \sum_{j=1}^n \frac{m_{j,1}}{1 - \rho_j} (1 - \rho_j) \rho_{[j,n]}
o(1) \\
\to & 
\frac{(k-1)!}{2^{k-1}} \lambda \nu^{k-1} \\
= & (k-1)! \lambda \left( \frac{\nu}{2} \right)^{k-1}.
\end{align*}
So for the sum \eqref{sumG} we have
\begin{eqnarray*}
\sum_{j=1}^n m_{j,1} \left( \overline G_{j+1,n}(x) - 1 \right)
& = & \sum_{j=1}^n m_{j,1} \sum_{k=1}^\infty
\frac{\overline G_{j+1,n}^{(k)}(1)}{k!} (x-1)^k \\
& = & \sum_{k=1}^\infty \frac{(x-1)^k}{k!}
\sum_{j=1}^n m_{j,1} \overline G_{j+1,n}^{(k)}(1) \\
& \to & \sum_{k=1}^\infty \frac{(x-1)^k}{k!} (k-1)!
\lambda \left( \frac{\nu}{2} \right)^{k-1} \\
& = & - \frac{2 \lambda}{\nu}
\log \left( 1 - (x-1)\frac{\nu}{2} \right),
\end{eqnarray*}
and Lemma 6 in \cite{GyIPV} makes the calculation rigorous.
Note that since $\overline G_{k,n}(x)$ is a polynomial the infinite sum
above is in fact finite.

We obtained that the limit generating function is
$$
\exp \left\{
- \frac{2 \lambda }{\nu} \log \left( 1 - (x-1)
\frac{\nu}{2} \right) \right\}
= \left( 1 - (x-1)\frac{\nu}{2} \right)^{- \frac{2 \lambda}{\nu}}
= \left( \frac{\frac{2}{2 + \nu}}{1 - \frac{\nu}{2 +  \nu} x}
\right)^{\frac{2 \lambda}{\nu}},
$$
which is the generating function of a negative binomial distribution with parameter
$r = 2 \lambda / \nu$ and $p = \nu/ (2 + \nu)$, as we stated.

\noindent \textbf{Step 2.}
Let us weaken the condition \eqref{G-assumption} on the generating functions.

From formula \eqref{comp-der} induction argument shows that
$$
\overline G_{j+1,n}^{(k)}(1)
= k! \sum_{s=2}^k  \sum_{l= j+1}^n \rho_{[j,l-1]} G_l^{(s)}(1)
\sum_{(\nu_1, \ldots, \nu_k)} \prod_{k=1}^m \frac{1}{\nu_m!}
\left( \frac{\overline G_{l+1,n}^{(m)}(1)}{m!} \right)^{\nu_m}.
$$
We claim that under assumption (iii) of the theorem \eqref{k-th der-asympt}
holds. By Step 1, this is true for $k=1,2$. Assume that the statement is
true until some $k \geq 2$. The previous case shows that we get the
asymptotic from the term corresponding to $s=2$, and we show that the terms
corresponding to  $s \geq 3$ are $o(\rho_{[j,n]})$. Using the induction
hypothesis for the term corresponding to $s \geq 3$ above we have
\begin{eqnarray*}
&& \sum_{l=j+1}^n \rho_{[j,l-1]} G_l^{(s)}(1)
\sum_{(\nu_1, \ldots, \nu_k)} \prod_{k=1}^n \frac{1}{\nu_m!}
\left( \frac{\overline G_{l+1,n}^{(m)}(1)}{m!} \right)^{\nu_m} \\
&& \leq \sum_{l=j+1}^n \rho_{[j,l-1]} G_l^{(s)}(1)
\sum_{(\nu_1, \ldots, \nu_k)} \prod_{k=1}^m \mathrm{const} \,
\rho_{[l,n]}^{\nu_m} \\
&& = \mathrm{const} \, \rho_{[j,n]} \sum_{l=j+1}^n
\frac{G_l^{(s)}(1)}{\rho_l ( 1 - \rho_l)} ( 1 - \rho_l) \rho_{[l,n]}^{s-1}  \\
&& = \rho_{[j,n]} o(1).
\end{eqnarray*}
We proved \eqref{k-th der-asympt}, and the statement follows as in the
previous case.

\noindent \textbf{Step 3.}
Finally, in case of general immigration distribution, according to the already proved part
what we have to show is
\begin{eqnarray*}
&&\left| \prod_{j=1}^n H_j \left( \overline G_{j+1,n}(x) \right) -
\prod_{j=1}^n
\left[ 1 + m_{j,1} \left( \overline G_{j+1,n}(x) - 1\right) \right] \right| \\
&& \hspace{70pt} \leq \sum_{j=1}^n \left| H_j \left( \overline G_{j+1,n}(x) \right)-
\left[ 1 + m_{j,1} \left( \overline G_{j+1,n}(x) - 1\right) \right] \right| \\
&& \hspace{70pt} \leq \sum_{j=1}^n \frac{m_{j,2}}{2} \left( \overline G_{j+1,n}(x) - 1\right)^2 \\
&& \hspace{70pt} \leq \frac{( x - 1)^2}{2} \sum_{j=1}^n m_{j,2} \rho_{[j,n]}^2
\to 0,
\end{eqnarray*}
where we used Lemma 5 and 6 in \cite{GyIPV}, and the assumption $m_{n,2}/(1 - \rho_n) \to 0$.
The proof is complete.
\end{proof}

\begin{proof}{Theorem \ref{th-gyors-konv}}
The generating function of the $n^{\textrm{th}}$ generation is 
$$
F_n(x) = \prod_{j=1}^n H_j( \overline G_{j+1,n}(x) ).
$$

First we show the existence of $\overline G_{j,\infty}(x)$.
Let us fix $j$, and investigate $\overline G_{j,n}(x)$ as $n \to \infty$.
Using the definition and the monotonicity of $\overline G_{j,n}$ and
that $G_n(x) \geq x$ for all $n$, we have
$$
\overline G_{j,n+1}(x) = \overline G_{j,n}(G_{n+1}(x))
\geq \overline G_{j,n}(x).
$$
Since $\overline G_{j,n}(x) \leq 1$, the limit
$$
\overline G_{j,\infty}(x) = \lim_{n \to \infty } \overline G_{j,n}(x)
$$
exists. Moreover $\overline G_{j+1, n} (x) \geq 1 - \rho_{[j,n]} ( 1- x)$
implies
$$
\overline G_{j+1, \infty} (x ) \geq 1 - \frac{\rho ( 1- x)}{\rho_{[1,j]}}.
$$
Therefore we can define the function $g$ in the theorem.

Using the inequality above, the convergence $m_{j,1} \to 0$,
the estimation $H_j(x) \geq 1-  m_{j,1} (1-x)$ and that for
$x$ small enough $1 - x \geq \e^{-(1 + \varepsilon) x}$ we obtain
\begin{eqnarray*}
H_j ( \overline G_{j+1,\infty}(x) )
& \geq &
1 -  m_{j,1} \left( 1 - \overline G_{j+1,\infty}(x) \right)
\geq 1 - m_{j,1} \frac{\rho ( 1- x)}{\rho_{[1,j]}} \\
& \geq & \exp \left\{ - (1 +  \varepsilon)
m_{j,1} \frac{\rho ( 1- x)}{\rho_{[1,j]}} \right\}.
\end{eqnarray*}
This easily implies that the function $g$ is continuous at 1,
$g(1) = 1$.

Next we show that $F_n(x) \to g(x)$, for all $ x \in [0,1]$. Introduce the notation
$$
g_n(x) = \prod_{j=1}^n H_j( \overline G_{j+1, \infty}(x)).
$$
Clearly $g_n(x) \to g(x)$, so we only have to show the convergence
$g_n(x) - F_n(x) \to 0$. We have
\begin{eqnarray*}
\left| g_n(x) - F_n(x) \right|
& \leq & \sum_{j=1}^n \left| H_j(\overline G_{j+1,\infty}(x))
- H_j(\overline G_{j+1,n}(x)) \right| \\
& \leq & \sum_{j=1}^n m_{j,1}
\left| \overline G_{j+1,n}(x) - \overline G_{j+1, \infty} (x) \right| \\
& \leq & \sum_{j=1}^{n_0}  m_{j,1}
\left| \overline G_{j+1,n}(x) - \overline G_{j+1, \infty} (x) \right| +
2 \sum_{j=n_0 +1}^\infty m_{j,1},
\end{eqnarray*}
where the first term goes to 0 for every fixed $n_0$, while the second one
can be arbitrary small by choosing $n_0$ large enough.
\end{proof}

\bigskip

\noindent \textbf{Acknowledgement}

I am grateful to Vincent Bansaye, who suggested to take a look at the linear
fractional case, which eventually lead to the final form of
Theorem \ref{th-nb}.

\end{document}